\documentclass[12pt]{article}

\usepackage{amsmath}  
\usepackage{amsfonts}  
\usepackage[english]{babel}  
\usepackage[T1]{fontenc}   
\usepackage{a4wide}   
   
\numberwithin{equation}{section}

\newtheorem{prop}{Proposition}[section]  
\newtheorem{corollary}{Corollary}[section]  
  
\newtheorem{lemma}{Lemma}[section]  
\newtheorem{Rem}{Remark}[section]

\newcommand{\CQFD}{ \hfill$\square$} 
  
\newcommand{\ind}{\mathbf{1}}  
  
\newcommand{\disp}{\displaystyle}  
\newcommand{\cvar}{\stackrel{var}{\longrightarrow}}  
\newcommand{\lto}{\longrightarrow}

\newcommand{\om}{\omega}  
  
\newcommand{\ee}{\mathbb{E}}

\newcommand{\retirer}[1]{$ $\newline  cf. version générale.}

\def\rit{\mathbb{R}}  
  
\def\nit{\mathbb{N}}

\def\pit{\mathbb{P}}  
\def\E{\mathop{\hbox{\rm I\kern-0.20em E}}\nolimits}

\def\med{\hskip 10pt}

\def\og{\leavevmode\raise.3ex  
     \hbox{$\scriptscriptstyle\langle\!\langle$~}}  
\def\fg{\leavevmode\raise.3ex  
     \hbox{~$\!\scriptscriptstyle\,\rangle\!\rangle$}~}

\allowdisplaybreaks

\def\bom{\bar \omega}
\date{ }  
  
\title{Regularity of the laws of shot noise series and of related processes
\footnote{This is a version of an original paper to be published in the Journal of Theoretical Probability
 which only differs from the published paper by typographical changes.}
}  

\author{
Jean-Christophe Breton
\\
\small{Laboratoire de Mathématiques et Applications}\\
\small{Universit\'e de La Rochelle, 17042 La Rochelle Cedex, France}\\
\small{Email: jcbreton@univ-lr.fr}\\
}

\begin{document}  

\maketitle 
\begin{abstract}
We investigate the regularity of shot noise series and of Poisson integrals. 
We give conditions for the absolute continuity of their law with respect to Lebesgue measure 
and for their continuity in total variation norm. In particular, the case of truncated series in adressed. 
Our method relies on a disintegration of the probability space based on a mere conditioning by the first jumps of the underlying Poisson process. 

\medskip \noindent
{\bf Keywords:} Absolute continuity, Convergence in variation, Disintegration of probability space, 
Poisson integrals, Shot noise series.  
\\
{\bf 2000 MR Subject Classification.} {\bf Primary: }60B10, {\bf Secondary:} 60F99, 60H05, 60H10.
\end{abstract}




\section*{Introduction}
The regularity of the laws of functionals of Lévy process has been intensively studied since the early eighties. 
The main idea has been to introduce on the path-space (supporting jumps) of the Lévy process some smooth structure to apply next a stochastic calculus of variation on this structure. 
Since the "Malliavin calculus" tackles successfully the regularity (existence of density, smoothness of density, etc) for Wiener functionals (see for instance \cite{Nualart}), 
one approach, originally by Bismut in \cite{Bismut}, has been to extend this Malliavin calculus to more general probability space, 
introducing in particular a chain-rule of differentiation and an integration-by-part in the path-space, see \cite{BGJ}. 
This setting applies in particular when the Lévy measure of the Lévy process has a density with respect to the Lebesgue measure. 
Another approach is based on perturbations introduced by Picard (see \cite{Picard}) consisting in adding points to random configurations driving the functionals. 
This line of work differs from Bismut's one both from the viewpoint of the technique used 
and of the conditions required, for instance the Lévy measure is not required to possess a density for Picard's approach. 
More in the spirit of Bismut's approach, another way to introduce a smooth structure in a probability space has been given by Davydov with the "stratification method". 
In this case, the smooth structure comes from a measurable group of admissible transformations of the space (the admissibility is referred to the law under study), see \cite{DL} and \cite{DLS}.  
Via a disintegration of the probability space into layers (the orbits of the group), this approach reduces the analysis of the law on path-space to finite-dimensional distributions.  
This method supports in particular mild differentiability properties for the functionals and applies in particular in a Poissonian setting. 
Recently, in a cycle of papers \cite{Kulik1}, \cite{Kulik2}, \cite{Kulik3}, \cite{Kulik4}, Kulik has developped a modification of the stratification method well designed for the study of functionals of Poisson point measure 
(absolute continuity and convergence in variation of the laws).  
He applies a so-called "time-stretching transformation" by changing time, stretching every infinitesimal segment of time $dt$ by $e^{h(t)}$ for some suitable function~$h$. 

As it appears in the time-stretching transformation, the method of Kulik relies on some perturbation of the jumping times of the process. 
Actually perturbation in the jumps of the Lévy process is a common feature of all those methods. 
For instance, in \cite{BGJ,DL, DLS} the values of the jumps are modified;
in \cite{Picard}, a new jump is added to the path of the process; 
while, as seen above, the jumping times are changed in \cite{Kulik1}--\cite{Kulik4}.  
This common feature implies that all these approaches require some kind of differentiability with respect to the perburbation used 
and consequently non-trivial regularity on the functionals or/and on the Lévy measure are required. 

Another technique to study the absolute continuity of the law of Lévy functionals which does not rely on differentiability behavior is given by Sato in the monograph \cite{Sato}. 
It is based on a decomposition of the underlying Lévy process in two independent components, each equipped with a truncation of the initial Lévy measures. 
The regularity  of the law of the functionals usually comes from the component equipped with the finite truncated Lévy measure which is a compound process whose law may be easily analysed. 
We refer for instance to Theorem 27.7 in \cite{Sato} and to \cite{PZ} for a recent application of such decomposition to the study of the law of Ornstein-Uhlenbeck processes.

In this paper we deal with specific Lévy functionals having some decomposition in series, the so-called "shot noise series", see \cite{Rosinski},
and we investigate the regularity of its law (by  regularity, we mean throughout either absolute continuity or convergence in variation). 
Our method is a version of the stratification method from \cite{DLS} and relies, like in \cite{PZ, Sato}, on a decomposition of the series in two parts. 
However in contrast to Sato's technique whose decomposition is based on the jump size, our decomposition is based on the jumping times: 
it is obtained by conditioning by some jumping time and yields a "regular part" (related to the first jumps and from which the regularity will be derived) 
and an "interfering part" (related to the last jumps).  
For instance in the most simple case, conditioning by (say) the second jump allows to derive the regularity of the law of the whole series from the behavior of the first summand, 
in a pretty elementary way. 
Similar type of techniques inspired from \cite{DLS} have been also recently used in \cite{Simon} to generalize \cite{PZ} for Ornstein-Uhlenbeck processes.
In contrast with the literature quoted above \cite{BGJ, Bismut, DL, DLS, Kulik1, Kulik2, Kulik3, Kulik4, Picard}, 
we do not apply any pertubation on the jump and 
our conditions are expressed directly in terms of some image measure of the kernel appearing in the shot noise series. 
Specifically, the differentiability of the kernel is not required in our results. 
Kulik's results for instance have a broader scope and they apply in particular for shot noise series (see Theorem 3.1 in \cite{Kulik1}, Theorems 2.1 and 2.2 in \cite{Kulik4}). 
But, as indicated above, it requires the differentiability of the kernel. 
However, of course even in our results, the differentiability is useful to state explicit and tractable sufficient conditions, see Lemma \ref{prop:ABP} and Remarks \ref{rem:P1a1}, \ref{rem:var1b} below.

The rest of the paper is organized as follows. 
In Section \ref{sec:shotnoise}, we introduce precisely the shot noise series, recalling conditions under which it is well defined. 
We give some modeling interpretations of those Lévy functionals and we exhibit its connections with Poisson integrals. 
In Section \ref{sec:AC}, we investigate the absolute continuity of the law of shot noise series. 
Section \ref{sec:var} is devoted to the continuity for the total variation norm of the laws of shot noise series with respect to their kernel.  
Note that in the whole paper, we deal also systematically with "truncated" shot noise series for which specific problems arise. 
Finally in Section \ref{sec:EDS}, we use similar techniques to study the regularity (in variation) of the law of simple Lévy-driven stochastic differential equation. 

In the whole sequel, we denote $\mu_1\ll\mu_2$ for the absolute continuity of a measure $\mu_1$ with respect to a measure $\mu_2$ 
and $\mu_1\asymp \mu_2$ for the equivalence of $\mu_1$ and $\mu_2$, {\it i.e.} $\mu_1\ll\mu_2$ and $\mu_2\ll\mu_1$. 
Implicitely, the absolute continuity of a measure $\mu$ refers to the Lebesgue measure $\lambda$. 
The restriction of the Lebesgue measure to some Borelian set $A$ is denoted by $\lambda_A$. 
The total variation of a signed measure $\mu$ is given by $\|\mu\|=\mu_+(A_\ast)+\mu_-(A_\ast^c)$ where
$\mu_+$, $\mu_-$ are the unique positive measures in the Hahn's decomposition of $\mu$ ($\mu=\mu_+-\mu_-$) 
and $A_\ast$ is such that $\mu_+(A_\ast)=\mu_+(\rit)$, $\mu_-(A_\ast^c)=0$. 
The corresponding convergence is denoted $\cvar$. 
Finally, we note ${\cal L}(X)$ for the law of a random variable $X$ or $X\sim\mu$ to indicate that $\mu$ is the law of $X$.


\section{Shot noise series and Poisson integrals}
\label{sec:shotnoise}
Let $(\Delta_i)_{i\geq 1}$ be a sequence of independent and identically distributed ({\it i.i.d.}) random vectors in $\rit^d$ with common law $\sigma$ without atom in $0$ 
and, independently, let $(T_i)_{i\geq 1}$ be a sequence of partial sums of independent exponential random variables with parameter $\alpha$.
Given a measurable function $h: \rit_+\times \rit^d \lto \rit$, we are interested in this paper in the regularity of the law of the so-called shot noise series
\begin{equation}
\label{eq:series}
I=\sum_{i=1}^{+\infty} h(T_i,\Delta_i).
\end{equation}
Namely, we give conditions on $h$ and on the law $\sigma$ for the absolute continuity 
of the law of $I$ in \eqref{eq:series} on $\rit$ or on $\rit\setminus\{0\}$.
We study also the continuity of the law for the total variation norm with respect to the kernel $h$.
Typically, our conditions are stated in terms of image measure $(\lambda\otimes \sigma)h^{-1}$ 
and we refer to \cite{ABP}, \cite{Davydov3} and \cite[Section 4]{DLS} for sufficient explicit conditions on $h$ and on $\sigma$ ensuring our conditions. 
\bigskip

The interest of such random series lies in their connections with Poisson integrals and in their interpretation as "shot noise series". 
Namely, the series \eqref{eq:series} can be seen as the cumulated effect of a repeated signal $\Delta_i$ under a "filter" $h$. 
More precisely, $h(T_i,\Delta_i)$ can be viewed as an effect at time zero of a shot $\Delta_i$ that happened $T_i$ time units ago
and the series \eqref{eq:series} represents the total cumulated shot noise at time zero.
In this setting, it is standard 
to suppose that $t\mapsto |h(t,x)|$ is non-increasing for all $x$, 
so that the magnitude of the effect decreases as the time elapsed from the moment of the shot increases. 
When $d=1$, a typical choice of filter is $h(t,x)=xg(t)$, with a fading effect as time goes on given by $g$ when it is decreasing. 
Such a case is related to Lévy-type integral (see below) and, in the sequel, we shall specify our results in this case. 
We refer to \cite{Rosinski} for a more precise interpretation of shot noise series.  
Shot noise series are also used to model internet backbone traffic. 
In such models, internet sessions are randomly opened according to Poisson arrivals $T_i$ and for a random duration $\Delta_i$. 
In this setting, the series $I$ represents an aggregate of packet streams generated by the random sessions, 
see \cite{BTIDO} for an analysis of flow-level information on an IP (Internet Protocol) backbone link, using Poisson shot noise series. 

\bigskip

We recall from Theorem 4.1 in \cite{Rosinski} necessary and sufficient conditions for the almost sure convergence of the series in \eqref{eq:series}, 
namely 
\begin{equation}
\label{eq:Lévy2}
\int_0^{+\infty}\int_{\rit^d\setminus\{0\}} (|h(t,x)|^2\wedge 1)dt\sigma(dx)<+\infty 
\end{equation}
and the existence of the following limit 
\begin{equation}
\label{eq:limit}
a(h):=\lim_{s\to+\infty} \int_0^s\int_{\rit^d\setminus\{0\}} h(t,x) \ind_{|h(t,x)|\leq 1} dt \sigma(dx).
\end{equation}
In particular, \eqref{eq:Lévy2}--\eqref{eq:limit} hold true and the series in \eqref{eq:series} is well defined when $h\in L^1(\lambda\otimes \sigma)$. 
In this case $\ee[|I|]\leq \int_0^{+\infty}\int_{\rit^d\setminus\{0\}} |h(t,x)|dt\sigma(dx)$. 
Moreover, recall that the law of the shot noise series in \eqref{eq:series} is infinitely divisible 
and its characteristic function is given by (see again \cite[Th. 4.1]{Rosinski}) 
\begin{equation} 
\label{eq:FCI}
\phi_I(u)=\exp\left(ia(h)u+\int_0^{+\infty}\int_{\rit^d\setminus\{0\}} \left(e^{iuh(t,x)}-1-iuh(t,x)\ind_{|h(t,x)|\leq 1}\right) dt\sigma(dx)\right).
\end{equation}

\bigskip
As indicated above, the shot noise series have also natural connection with Poisson integrals (more generally with Lévy-type integrals). 
We refer to \cite{Sato} for a precise description of Poisson and Lévy integrals. 
Let us fix some notations by describing the connections between shot noise series and Poisson integrals.  
For instance, for a Poisson measure $N$ on $\rit^d$ with control measure $\nu$, let $Z(t)=\int_0^t\int_{|x|\geq 1} xN(ds,dx)$ be the multidimensional Poisson point process associated to $N$, 
{\it i.e.}  $Z(t)=(Z_1(t), \dots, Z_d(t))$ with $Z_i(t)=\int_0^t\int_{|x|\geq 1} x_iN(ds,dx)$. 
Define $(T_i)_{i\geq 0}$ to be the sequence of ordered jumping times of $Z$ in $B_1^c=\{x\in\rit^d: |x|\geq 1\}$, 
that is $T_0=0$ and $T_i=\inf\{t> T_{i-1} : | \Delta Z_t|>1\}$ for $i\geq 1$. 
As a sum of $i$ independent exponential random variables with parameter $1/\nu(B_1^c)$, 
$T_i$ has a law equivalent to $\lambda_{\rit_+}$. 
Moreover, the jumps  $\Delta_i:=\Delta Z_{T_i}$ of $Z$
form an {\it i.i.d.} sequence of random variables with common law $\nu(\cdot \cap B_1^c)/\nu(B_1^c)$ 
and which is independent of the jumping times $(T_i)_{i\geq 0}$.
In this case, taking 
\begin{equation}
\label{eq:choix}
h(t,x)=f(t,x)\ind_{|x|\geq 1},\quad
\alpha=1/\nu(B_1^c), \quad 
\sigma=\nu(\cdot \cap B_1^c)/\nu(B_1^c),
\end{equation}
the series in \eqref{eq:series} becomes a representation for 
\begin{equation}
\label{eq:IP}
\int_0^{+\infty} \int_{|x|\geq 1}f(t,x) N(dt,dx). 
\end{equation}

In the particular case $f(t,x)=xg(t)$ when $d=1$, note that the Poisson integral in \eqref{eq:IP} is the non-compensated part of the Lévy-type integral 
\begin{equation}
\label{eq:gY}
\int_0^{+\infty} g(s) dY_s=
\int_0^{+\infty}\int_{|x|<1} xg(s) \widetilde N (ds,dx)+ \int_0^{+\infty}\int_{|x|\geq 1} xg(s) N(ds,dx) 
\end{equation}
where $Y_t$ is the Lévy process $\int_0^t\int_{|x|<1} x \widetilde N (ds,dx)+ \int_0^t\int_{|x|\geq 1} x N(ds,dx)$ 
and where $\widetilde N =N-\nu$ stands for the compensated Poisson measure.  
The shot noise series \eqref{eq:series} thus gives an insight into the law of the Lévy-integral. 
For instance, since the integrals in the right-hand side of \eqref{eq:gY} are independent,
the absolute continuity with respect to the Lebesgue measure of the law of $I$ ($=\int_0^{+\infty}\int_{|x|\geq 1} xg(s) N(ds,dx)$ in this case)
ensures that of $\int_0^{+\infty} g(s) dY_s$.
Of course, such a remark is useful only if there is no Gaussian part in the Lévy process $Y$. 
\\

When $h(s,x)=f(s,x)\ind_{[0,t]\times B_1^c}(s,x)$ or $h(s,x)=f(s,x)\ind_{[0,t]\times (\rit^d\setminus\{0\})}(s,x)$ with $\nu(\rit\setminus\{0\})<+\infty$, 
the shot noise series represent Poisson integrals computed on $[0,t]$. 
In this case, there is almost surely (a.s.) a finite number of summand in \eqref{eq:series} and the series is actually truncated 
\begin{equation}
\label{eq:series2}
I(t)=\sum_{\begin{subarray}{c} i\geq 1\\ T_i\leq t\end{subarray}} h(T_i,\Delta_i).
\end{equation}

We shall see that specific problems arise for truncated series (for instance, there is an obvious atom in the law at zero) and this case requires a deeper analysis. 
Note that, according to our interpretation in terms of cumulated shots, $I(t)$ represents the total cumulated effect of random shots than happened less than $t$ units of time before zero.  

Note that in this Poissonian interpretation, the shot noise series appears as functionals of Poisson point process.  
Moreover, following this interpretation, in the sequel, $T_i$'s will be referred to as the jumping times and the $\Delta_i$'s as the (corresponding) jumps. 
\bigskip


In our study of shot noise series $I$, a crucial step in the sequel is to condition for some $p>1$ by the jumping time $T_{p+1}$
and to derive the regularity of the series $I$ from the behavior of its first summands.
Typically, 
$p$ will be set to $1$. 
To do so, we shall consider (without restriction) that the probability space $(\Omega, {\cal F}, \pit)$ is disintegrated as follows:
\begin{equation}
\label{eq:desintegre}
\left(\bar \Omega_{p+1}\times[0, T_{p+1}(\bom)]^p\times (\rit^d\setminus\{0\}), {\cal F}_{p+1}^*\times {\cal B}([0, T_{p+1}(\bom)]^p)\times {\cal B}(\rit\setminus\{0\}), \bar \pit_{p+1}\otimes \bar\lambda_{[0, T_{p+1}(\bom)]}^{stat}\otimes \sigma\right)
\end{equation}
where we note ${\cal F}_{p+1}^*=\sigma(T_i, \Delta_i : i\geq p+1)$ and 
$\bar\lambda_{[0, T_{p+1}(\bom)]}^{stat}$ is the $p$-dimensional distribution of the order statistics on $[0,T_{p+1}(\bom)]$, 
{\it i.e.} the measure with density $(p!/T_{p+1}(\bom))\ind_{0\leq t_1\leq\dots\leq t_p\leq T_{p+1}(\bom)}$.
For any random functional $F$ on this space and every Borelian set $A$, we have:
\begin{eqnarray*}
&&\pit(F\in A)\\
&=&\pit(F(\bom,T_1, \Delta_1, \dots, T_p, \Delta_p) \in A)\\
&=&\bar \ee_{p+1}\left[\frac 1{T_{p+1}(\bom)^p}\int_{(\rit^d\setminus\{0\})^p}\int_{0\leq t_1\leq\dots\leq t_p\leq T_{p+1}(\bom)} 
\ind_{\{F(\bom, t_1,x_1, \dots t_p, x_p)\in A \}} \prod_{i=1}^pdt_i \sigma(dx_i)\right]
\end{eqnarray*}
where $\bar \ee_{p+1}$ stands for the expectation with respect to $\bar\pit_{p+1}$.
In the sequel, $h$ stands for a typical kernel in the shot noise series \eqref{eq:series} we consider. 
We use the notation $f$ when the kernel is, in a way, truncated, for instance 
$h(s,x)=f(s,x)\ind_{B_1^c}(x)$, $h(s,x)=f(s,x)\ind_{[0,t]}(s)$  or $h(s,x)=f(s,x)\ind_{[0,t]\times B_1^c}(s,x)$
and we use the notation $g$ when the kernel $h(t,x)=xg(t)$ is related to Lévy-type integral.


\section{Absolute continuity of the law of shot noise series}
\label{sec:AC}

In this section, we address the problem of the absolute continuity of the law of shot noise series. 
Since $I$ is infinitely divisible distributed without Gaussian part, it is known that it is not easy to find conditions ensuring absolute continuity. 
Nevertheless, such conditions are given for instance by Sato in \cite{Sato} for general infinite divisible law  
and by Kulik in \cite[Th. 2.1]{Kulik4} for functionals of Poisson point process. 
In particular, in our setting, with $\tilde\nu(dx)=(x^2\wedge 1)(\lambda\otimes\sigma)h^{-1}(dx)$, Sato shows in \cite[Th. 27.7]{Sato} that the law of $I$ in \eqref{eq:series} is absolutely continuous 
if, for some $p>0$, $\tilde\nu^{\ast p}\ll\lambda$. 
However, it is worth noting that several interesting shot noise series such as truncated one in \eqref{eq:series2} have an obvious atom in zero 
(the sum in \eqref{eq:series2} is empty on the non-negligible event $\{T_1> t\}$ so that $I(t)=0$). 
In this case, the general results are no longer applicable. 
In this section, we are interested in the absolute continuity on $\rit\setminus\{0\}$ of the law of shot noise series, see Proposition \ref{prop:P1b}. 
Our argument is simple and relies on the behavior of the first summands in the series by conditioning.
Since such truncated series are related to Poisson integrals on $[0,t]$ appearing for instance in some Lévy-type integrals, 
Proposition \ref{prop:P1b} gives an insight in the regularity of the law of Lévy-type integrals, see Corollary \ref{cor:P5}. 
\\

\noindent
First, in order to illustrate our technique in a simple context, let us recover a condition close to that of \cite[Th. 27.7]{Sato}, stated above.
Since for each $n$, the law of $T_n$ is equivalent to $\lambda_{\rit_+}$, the condition 
$(\lambda_{\rit_+}\otimes \sigma)h^{-1}\ll \lambda$
ensures the absolute continuity of the law of each summands $h(T_i,\Delta_i)$ in \eqref{eq:series}. 
Since the summands are not independent, this does not imply directly the absolute continuity of the law of $I$.
However because of the conditional independence of the summands, this appears to be sufficient. 
Actually, it is enough to suppose that, for some $p\in\nit\setminus\{0\}$,
\begin{equation}
\label{eq:H1p}
((\lambda_{\rit_+}\otimes \sigma)h^{-1})^{\ast p}\ll \lambda
\end{equation}
that is, roughly speaking, the absolute continuity of the law of the sum of any of the first terms ensures the absolute continuity of the total sum.
This is a strictly weaker condition than $(\lambda_{\rit_+}\otimes \sigma)h^{-1}\ll \lambda$ (see Remark 27.2 in \cite{Sato}) which recovers, up to a density, the condition in \cite[Th. 27.7]{Sato}.
Using the disintegrated probability space $\bar \Omega_{p+1}\times[0, T_{p+1}(\bom)]^p\times (\rit^d\setminus\{0\})^p$ defined in \eqref{eq:desintegre}, 
the proof of the absolute continuity of $\mu:={\cal L}(I)$ under \eqref{eq:H1p} comes from the absolute continuity of the conditional measures $\mu_{\bom}={\cal L}(I|{\cal F}_{p+1}^*)(\bom)$ in 
the decomposition $\mu=\int_{\Bar\Omega_{p+1}} \mu_{\bom}\: \bar\pit_{p+1}(d\bom)$ of $\mu$ as the mixture of the $\mu_{\bom}$. 
Indeed, with $\Sigma_{p+1}(\bom):=\sum_{i\geq p+1} h(T_i, \Delta_i)$ and 
$(T_1', \dots, T_p')\stackrel{{\cal L}}{=}{\cal L}(T_1, \dots, T_p|{\cal F}_{p+1}^*)(\bom)$ the uniform order statistics on $[0,T_{p+1}(\bom)]$, 
we have 
\begin{equation}
\label{eq:B1}
\mu_{\bom}={\cal L}\left(\sum_{i=1}^ph(T_i', \Delta_i)+\Sigma_{p+1}(\bom)\right).
\end{equation} 
Since $\Sigma_{p+1}(\bom)$ is known when $\bom\in\bar\Omega_{p+1}$ is given, the absolute continuity of $\mu_{\bom}$ is equivalent to that of 
${\cal L}(\sum_{i=1}^ph(T_i', \Delta_i))=(\bar\lambda_{[0, T_{p+1}(\bom)]}\otimes\sigma)h^{-1})^{\ast p}$.
This ensures $\mu_{\bom}\ll\lambda$ for all $\bom\in \bar \Omega_{p+1}$ and $\mu\ll\lambda$. 
 
\begin{Rem}
\label{rem:P1a1} 
{\rm In the sequel, we shall express our conditions in terms of regularity of image measure $(\lambda_{\rit_+}\otimes \sigma)h^{-1}$ like in \eqref{eq:H1p}
We give here some explicit sufficient conditions to ensure such a regularity.
Let $h_t$ and $h_x$ be the sections of $h$ defined by $h(t,x)=h_t(x)=h_x(t)$ 
and observe that
$$
(\lambda_{\rit_+}\otimes\sigma)h^{-1}(A)
=\int_{\rit^d} \lambda_{\rit_+}h_x^{-1}(A) \sigma(dx)
=\int_{\rit_+} \sigma h_t^{-1}(A) dt. 
$$
The condition \eqref{eq:H1p} is satisfied for $p=1$ if, for $\sigma$-almost all $x\in\rit^d$, $\lambda_{\rit_+}h_x^{-1}\ll\lambda$
or if, for almost all $t\in\rit_+$, $\sigma h_t^{-1}\ll\lambda$.
This is satisfied in particular if for $\sigma$-almost all $x\in\rit^d$, $h_x$ is differentiable a.e. with 
$\partial_t h(t,x)\not=0$ for almost all $t\in\rit$,
resp. if $\sigma\ll\lambda^d$ and for almost all $t\in\rit_+$, $h_t$ is differentiable a.e. with $\det D_x h(t,x)\not=0$ for $\sigma$-almost all $x\in\rit_+$, 
see Theorem 4.3 in \cite{DLS}. 
Here $\partial_t h(t,x)$ stands for the partial derivative of $h$ with respect to $t\in\rit_+$ 
and  $D_x h(t,x)$ for the matrix of first order partial derivatives of $h(t,x)$ with respect to $x\in\rit^d$. 

Such differentiability-type conditions along time for $h(t,x)$ are also required in the "time-stretching" method used by Kulik, see in particular Theorem 2.1 in \cite{Kulik4}. 
}
\end{Rem}

\bigskip
We investigate now the law of the shot noise series beside a possible atom in $0$. 
The typical interesting case is that of a truncated shot noise series $I(t)$ as defined in \eqref{eq:series2}. 
\begin{prop}
\label{prop:P1b}
Suppose that 
\begin{equation}
\label{eq:H2}
((\lambda_{\rit_+}\otimes\sigma)h^{-1})_{|\rit\setminus\{0\}} \ll \lambda, 
\end{equation} 
then the law of $I$ in \eqref{eq:series} is absolutely continuous with respect to $\lambda$ on $\rit\setminus\{0\}$, with possibly an atom at zero.
\end{prop}

\begin{Proof}
Consider the sequence of random variables $N_i(\om)=\min\{j\geq N_{i-1} : h(T_j,\Delta_j)\not=0\}$, with the convention $\min\{\emptyset\}=+\infty$.
This defines stopping times for the $\sigma$-algebrae $\sigma(T_i-T_{i-1}, \Delta_i : i\leq k)$, $k\in\nit\setminus\{0\}$.
Let $A_0=\{N_1=+\infty\}$, $A_1=\{N_1<N_2=+\infty\}$, $A_2=\{N_2<+\infty\}$. 
The law $\mu$ of $I$ splits as $\mu= \pit(A_0)\mu_0+\pit(A_1)\mu_1+\pit(A_2)\mu_2$
where $\mu_k=\pit(I\in \cdot |A_k)$ for $0\leq k\leq 2$.

We have $\mu_0=\delta_0$ since on $A_0$, the sum defining $I$ in \eqref{eq:series} is empty. 
On $A_1$, the sum defining $I$ in \eqref{eq:series} reduces to $h(T_{N_1}, \Delta_{N_1})$, 
and for any Borelian set $A$ with $\lambda(A)=0$, we have: 
\begin{eqnarray}
\nonumber
\mu_1(A)&=&\pit(h(T_{N_1}, \Delta_{N_1})\in A\:|\: A_1)\\
\nonumber
&=&\sum_{i_1=1}^{+\infty} \frac{\pit(h(T_{i_1}, \Delta_{i_1})\in A\setminus\{0\}, h(T_k, \Delta_k)=0 \ \forall k\neq i_1)}{\pit(A_1)}\\
\label{eq:A1}
&\leq &\sum_{i_1=1}^{+\infty} \frac{\pit(h(T_{i_1}, \Delta_{i_1})\in A\setminus\{0\})}{\pit(A_1)}.
\end{eqnarray}
But since ${\cal L}(h(T_{i_1}, \Delta_{i_1})) \asymp (\lambda_{\rit_+}\otimes\sigma)h^{-1}$, 
the condition \eqref{eq:H2} entails  $\mu_1(A)=0$, {\it i.e.} $\mu_1\ll \lambda$. 

\noindent
Next on $A_2$, first we condition by $N_2=i_2$ 
$$
\mu_2=\sum_{i_2=2}^{+\infty}  \frac{\pit(N_2=i_2)}{\pit(A_2)}\mu_{i_2}
$$
where $\mu_{i_2}={\cal L}(I|N_2=i_2)$.
Next, we condition by $T_{i_2}$: 
$$
\mu_{i_2}=\int \tilde \mu_{i_2}\ d\pit_{T_{i_2}} 
$$
where the measure $\tilde \mu_{i_2}$ is the law of 
$$
\sum_{k=1}^{i_2-1} h(T_k,\Delta_k)+\sum_{k=i_2}^{+\infty} h(T_{i_2}+T_k', \Delta_k)
=h(T_{N_1}',\Delta_{N_1}')+\sum_{k=i_2}^{+\infty} h(T_{i_2}+T_k', \Delta_k)
$$
where ${\cal L}(T_{N_1}',\Delta_{N_1})={\cal L}((T_{N_1},\Delta_{N_1})|N_2=i_2,T_{i_2})$ 
and for $k\geq i_2$, $T_k'=T_k-T_{i_2}$ is independent of $\sigma(T_k,\Delta_k : k\leq i_2-1)$.
Thus $\tilde \mu_{i_2}$ is absolutely continuous with respect to $\lambda$ if the conditonal law ${\cal L}(h(T_{N_1},\Delta_{N_1})|N_2=i_2, T_{i_2})$ is. 
But for a Borelian set $A$ with $\lambda(A)=0$, we have
\begin{eqnarray}
\nonumber
&&\pit(h(T_{N_1},\Delta_{N_1})\in A |N_2=i_2, T_{i_2})\\
\nonumber
&=&\sum_{i_1=1}^{i_2-1}\pit(h(T_{i_1},\Delta_{i_1})\in A, N_1=i_1 |N_2=i_2, T_{n_2})\\
\nonumber
&=&\sum_{i_1=1}^{i_2-1}\pit(h(T_{i_1},\Delta_{i_1})\in A\setminus \{0\}, h(T_k, \Delta_k)=0\ \forall i_1\neq k<i_2|N_2=i_2, T_{i_2})\\
\nonumber
&\leq&\sum_{i_1=1}^{i_2-1}\pit(h(T_{i_1},\Delta_{i_1})\in A\setminus \{0\} |N_2=i_2, T_{i_2})\\
\nonumber
&=&\sum_{i_1=1}^{i_2-1}\frac{\pit(h(T_{i_1},\Delta_{i_1})\in A\setminus \{0\}, N_2=i_2|T_{i_2})}{\pit(N_2=i_2|T_{i_2})}\\
\label{eq:A2}
&\leq&\sum_{i_1=1}^{i_2-1}\frac{\pit(h(T_{i_1}',\Delta_{i_1})\in A\setminus \{0\} )}{\pit(N_2=i_2|T_{i_2})}
\end{eqnarray}
where here ${\cal L}(T_{i_1}')={\cal L}(T_{i_1}|T_{i_2})$ is the $i_1$-th uniform order statistics on $[0,T_{i_2}]$. 
But since $h(T_{i_1}',\Delta_{i_1})\simeq (\bar\lambda_{[0,T_{i_2}(\bom)]}\otimes\sigma)h^{-1}$, 
condition \eqref{eq:H2} entails 
$\pit(h(T_{N_1},\Delta_{N_1})\in A |N_2=i_2, T_{i_2})=0$ 
and finally $\mu_2\ll \lambda$. 
We conclude the proof of Proposition \ref{prop:P1b} gathering all the intermediate results. 
\CQFD
\end{Proof}


\bigskip
In case $h(t,x)=xg(t)$ when $d=1$, the previous results can be specialized for the non-compensated part in Lévy-type integrals \eqref{eq:gY}. 
Because of the special structure of $h$, conditions for absolute continuity can be proposed independently of the measure $\sigma$. 
\begin{corollary}
\label{cor:P5}
Let $g:\rit^+\to\rit$ and consider the shot noise series $I$ in \eqref{eq:series} but with $h(t,x)=xg(t)$. 
\begin{enumerate}
\item The law of $I$ is absolutely continuous with respect to $\lambda$ on $\rit$ if: 
\begin{equation}
\label{eq:H1b}
\lambda_{\rit_+} g^{-1} \ll \lambda. 
\end{equation} 

\item The law of $I$ is absolutely continuous with respect to $\lambda$ on $\rit\setminus\{0\}$ if:
\begin{equation*}
\label{eq:H2b}
(\lambda_{\rit_+} g^{-1})_{|\rit\setminus\{0\}} \ll \lambda. 
\end{equation*} 
\end{enumerate}
\end{corollary}
\begin{Rem}
\label{rem:b}
{\rm 
\begin{itemize}
\item In order to state an "easy-to-check" condition in 1) as in \eqref{eq:H1b}, we do not give a strict counterpart of 
\eqref{eq:H1p} in terms of convolution. 
\item
Like in Remark \ref{rem:P1a1}, if $g$ is differentiable a.e. with $g'(t)\not=0$ a.e., 
condition \eqref{eq:H1b} is satisfied, see Theorem 4.2 in \cite{DLS}. 
\item 
Like for Proposition \ref{prop:P1b}, the second statement is interesting in particular for truncated shot noise series $I(t)$ for which we know there is an atom in zero. 
\end{itemize}
}\end{Rem}
\begin{Proof}
1) Follow the same lines as in the simple argument before Proposition \ref{prop:P1b} with $p=1$. 
Plugging $h(T_1',\Delta_1)=\Delta_1 g(T_1')$ in \eqref{eq:B1}, 
the proof reduces to the absolute continuity of the law of $\Delta_1 g(T_1')$.
But since $g(T_1')$ is absolutely continuous under \eqref{eq:H1b}, the conclusion comes from the following remark: 
 a product $XY$ of independent random variables $X, Y$ has an absolutely continuous law whenever $X$ has one and $Y$ does not have atom in zero ($\pit(Y=0)=0$). 
\\

2) We follow similarly the same lines as in the proof of Proposition \ref{prop:P1b}. 
Plugging $h(T_{i_1},\Delta_{i_1})=\Delta_{i_1} g(T_{i_1})$ in \eqref{eq:A1}
and $h(T_{i_1}',\Delta_{i_1})=\Delta_{i_1} g(T_{i_1}')$ in \eqref{eq:A2}, 
the proof reduces to the absolute continuity of the law of $\Delta_{i_1} g(T_{i_1}')$.
We conclude like in 1) with a similar remark: 
a product $XY$ of independent random variables $X, Y$ has also an absolutely continuous law on $\rit\setminus\{0\}$ law 
when ${\cal L}(X)_{|\rit\setminus\{0\}}\ll \lambda$ and $\pit(Y=0)=0$. 
\CQFD
\end{Proof}


\section{Convergence in variation of the law of shot noise series}
\label{sec:var}

In this section, we study further the law of shot noise series
by investigating its behavior for the total variation norm with respect to the "filter" $h$. 
Such convergences have been recently given by Kulik for general functionals of Poisson point process under some kind of differentiability of the functionals, see Theorem 2.2 in \cite{Kulik4}. 
Here, we propose in an elementary fashion convergence in variation under conditions expressed in terms of image measure by the filter $h$. 
Explicit sufficient conditions are available when the filters are smooth functions using Lemma \ref{prop:ABP} below from \cite{ABP}. 
See also Remark \ref{rem:var1b} for a comparison with \cite{Kulik4}.  
Note that when the laws of the shot noise series have densities, 
the convergence in variation is equivalent to the convergence in $L^1(\rit)$ of the densities.
In the sequel, we deal in Section \ref{sec:varI} with series $I$ in \eqref{eq:series}, 
typical examples are Poisson integrals \eqref{eq:IP} on $\rit_+\times B_1^c$ when \eqref{eq:choix} holds true. 
The case of truncated series $I(t)$ in \eqref{eq:series2} is more difficult since the laws have an atom in zero 
and this case is addressed in Section~\ref{sec:varI(t)}. 
We begin with some useful results on convergence in variation in Section~\ref{sec:varStaff}.


\subsection{On convergence in variation}
\label{sec:varStaff}
We shall use the following elementary results. 
For the sake of self containess, we include their proofs:
\begin{lemma}
\label{lemme:1}
Let $X, Y, Z$ be random variables such that $Z$ is independent of $(X,Y)$ and $\pit(Z=0)=0$. 
Then $\|{\cal L}(XZ)-{\cal L}(YZ)\|\leq \|{\cal L}(X)-{\cal L}(Y)\|$. 
\end{lemma}

\begin{Proof} 
Let $\pit_X$ stand for the law of the random variable $X$. 
Using the Hahn's decomposition of $\pit_{XZ}-\pit_{YZ}$ into $\pit_{XZ}-\pit_{YZ}=\mu_+-\mu_-$ 
and of $\pit_X-\pit_Y$ into $\pit_X-\pit_Y=\nu_+-\nu_-$,
we have 
$$
\|\pit_{XZ}-\pit_{YZ}\|=\mu_+(A_\ast)+\mu_-(A_\ast^c), \quad 
\|\pit_X-\pit_Y\|=\nu_+(B_\ast)+\nu_-(B_\ast^c)
$$
where $\mu_+(A_\ast)=\mu_+(\rit)$, $\mu_-(A_\ast^c)=0$ and 
$\nu_+(B_\ast)=\nu_+(\rit)$, $\nu_-(B_\ast^c)=0$.
Note that, for any Borelian set $A$,  
$
-\nu_-(B_\ast^c)\leq-\nu_-(A)\leq \pit_X(A)-\pit_Y(A)\leq \nu_+(A)\leq \nu_+(B_\ast)
$
and $\pit(XZ\in A)=\ee_{Z}[\pit_{X}(A_Z)]$ where $\ee_Z$ is the expectation with respect to $Z$ and for any $z\not=0$, $A_z=\{a/z : a\in A\}$ 
(note that $A_Z$ is almost surely well defined since $\pit(Z=0)=0$). 
We have  
\begin{eqnarray*}
\mu_+(A_\ast)&=&\pit_{XZ}(A_\ast)-\pit_{YZ}(A_\ast)=\ee_Z[\pit_X((A_\ast)_Z)-\pit_Y((A_\ast)_Z)]\leq \nu_+(B_\ast)\\
-\mu_-(A_\ast^c)&=&\pit_{XZ}(A_\ast^c)-\pit_{YZ}(A_\ast^c)=\ee_Z[\pit_X((A_\ast^c)_Z)-\pit_Y((A_\ast^c)_Z)]\geq -\nu_-(B_\ast^c)
\end{eqnarray*}
from which the conclusion of the lemma derives.
\CQFD
\end{Proof} 

\begin{lemma}
\label{lemme:convol}
Let $\mu$, $\nu$ be two finite (signed) measures. For any $p\in\nit\setminus\{0\}$, we have 
$$
\|\mu^{\ast p}-\nu^{\ast p}\|
\leq p\max(\|\mu\|,\|\nu\|)^{p-1}\|\mu-\nu\|.
$$
In particular, if for probability measures $\mu_n\cvar \mu$ when $n\to+\infty$, 
then $\mu_n^{\ast p}\cvar \mu^{\ast p}$ for any $p\geq 1$.
\end{lemma}
\begin{Proof}
The statement is obvious for $p=1$. 
Suppose it holds true for $p\geq 1$, we have 
\begin{eqnarray*}
\|\mu^{\ast (p+1)}-\nu^{\ast (p+1)}\|
&\leq& \|\mu\ast(\mu^{\ast p}-\nu^{\ast p})\|+\|(\mu-\nu)\ast\nu^{\ast p}\|\\
&\leq& \|\mu\|\|\mu^{\ast p}-\nu^{\ast p}\|+\|\nu^{\ast p}\|\|\mu-\nu\|\\
&\leq& p\|\mu\|\max(\|\mu\|,\|\nu\|)^{p-1}\|\mu-\nu\|+\|\nu\|^p\|\mu-\nu\|\\
&\leq& (p+1)\max(\|\mu\|,\|\nu\|)^p\|\mu-\nu\|.
\end{eqnarray*}
This proves the statement by induction. 
\CQFD
\end{Proof}

\bigskip 
Our main conditions for convergence of variation of the law of shot noise series are expressed
in terms of image measures (see \eqref{eq:cvar-main}, \eqref{eq:H3b}, \eqref{eq:H3bis}, \eqref{eq:H5} below).  
In order to give more explicit sufficient conditions in Remarks \ref{rem:var1}, \ref{rem:varg1}, \ref{rem:var2}, \ref{rem:varg2}, 
we shall use the following result. 
This is an adaptation, specially designed to our setting, of Theorem 3.1 and Corollary 3.2 in \cite{ABP} 
(where we set therein $X=\rit^{p+q}$, $Y=\rit^p$, $Z=\rit^q$ and $\mu=\lambda^{p+q}$, 
$\mu_0=\lambda^p\otimes \sigma$). 
We denote $\partial_{y_i}F(y,z)$ for the derivative of $F(y,z)$ with respect to the $i$-th coordinate of $y$, and the same for 
$\partial_{z_j}F(y,z)$.  
\begin{lemma}
\label{prop:ABP}
Let $F_n, F:\rit^{p+q}\to\rit$ such that for almost all $y$, 
$F_n(y,z)$ and $F(y,z)$ are absolutely continuous as functions of $z_i$, $1\leq i\leq q$. 
Suppose $F_n\to F$ and $\partial_{z_i} F_n(y,z)\to\partial_{z_i} F(y,z)$ in $L^r(\rit^{p+q},\lambda^{p+q})$ for $r\geq p+q$ when $n\to +\infty$.
Suppose moreover that for almost all $(y,z)$, $\det(\partial_{z_i} F(y,z),\partial_{z_j} F(y,z))_{i,j=1}^q\neq 0$
Then $\lambda^{p+q} F_n^{-1} \cvar \lambda^{p+q} F^{-1}$ when $n\to+\infty$. 
Moreover if $\sigma\ll \lambda^q$, we have also 
$(\lambda^p\otimes\sigma)F_n^{-1}\cvar (\lambda^p\otimes\sigma)F^{-1}$.  
\end{lemma}
We shall use also Lemma \ref{prop:ABP} in a one-dimensional setting. 
In this case, Lemma \ref{prop:ABP} reduces to the original result due to Davydov in \cite{Davydov3}. 
\begin{lemma}
\label{lemme:Davydov3}
Suppose that the functions $f_n$ and $f$ are absolutely continuous on $[\alpha,\beta]$ such that
$\lim_{n\to+\infty} f_n(\alpha)= f(\alpha)$,
$\lim_{n\to+\infty}\| f_n'-f'\|_{L^1([\alpha,\beta])}= 0$ and
$f'(x)\not=0$ for almost all $x\in[\alpha,\beta]$.
Then $\lambda_{[\alpha,\beta]} f_n^{-1}\cvar \lambda_{[\alpha,\beta]} f^{-1}$ , $n\to+\infty$. 
\end{lemma}

\subsection{Convergence of shot noise series}
\label{sec:varI}
In this section, we give conditions for the continuity in total variation norm of the law of shot noise series $I$ as in \eqref{eq:series} with respect to $h$. 
\begin{prop}
\label{prop:P4}
Let $h_n, h$ satisfy \eqref{eq:Lévy2}--\eqref{eq:limit}.
Suppose 
\begin{enumerate}
\item $\lim_{n\to+\infty}a(h_n)=a(h)$; 
\item $\lim_{n\to+\infty} \int_0^{+\infty}\int_{\rit^d\setminus\{0\}} (|h_n(t,x)-h(t,x)|^2\wedge 1) dt\sigma(dx)=0$; 
\item for some $p\in \nit\setminus\{0\}$, $((\lambda_{\rit_+}\otimes \sigma)h^{-1})^{\ast p}\ll\lambda$
 and for all $t>0$:
\begin{equation} 
\label{eq:cvar-main}
((\lambda_{[0,t]}\otimes \sigma)h_n^{-1})^{\ast p}\cvar ((\lambda_{[0,t]}\otimes \sigma) h^{-1})^{\ast p}.
\end{equation}
\end{enumerate}
Then $\disp {\cal L}(I_n)\cvar {\cal L}(I)$ when $n\to+\infty$.
\end{prop}
\begin{Proof}
From the main condition \eqref{eq:cvar-main}, we have 
\begin{equation}
\label{eq:mu1}
\tilde\mu_{n,p}:={\cal L}\left(\sum_{i=1}^ph_n(U,\Delta_i)\right)\cvar \tilde\mu_p:={\cal L}\left(\sum_{i=1}^ph(U,\Delta_i)\right)
\end{equation}
for any independent random variables $U_i$, $1\leq i\leq p$, uniform on $[0,t]$ and independent of $(\Delta_i)_{i\geq 1}$. 
We disintegrate the probability space as in \eqref{eq:desintegre}
and the total variation of measures rewrites:
\begin{equation}
\label{eq:desint1} 
\|\mu_n-\mu\|=\int_{\bar\Omega_{p+1}} \|\mu_{n,\bom}-\mu_{\bom}\|\bar\pit_{p+1}(d\bar\omega).
\end{equation}
Let $\tau_{n,\bom}(x)=x+\Sigma_{p+1}^n(\bom)$, resp. $\tau_{\bom}(x)=x+\Sigma_{p+1}(\bom)$, be the translation of $\Sigma_{p+1}^n(\bom)$, resp. of $\Sigma_{p+1}(\bom)$. 
The measure $\mu_{n,\bom}$ is the law of $\sum_{i=1}^pf_n(T_i', \Delta_i)+\Sigma_{p+1}^n$ 
where $(T_1', \dots, T_p')\stackrel{{\cal L}}{=}{\cal L}(T_1, \dots, T_p |T_{p+1}(\bom))$ is the uniform order statistic on $[0,T_{p+1}(\bom)]$.
It rewrites $\mu_{n,\bom}=\tilde \mu_{n,p} \tau_{n,\bom}^{-1}$ and we have  
\begin{eqnarray}
\nonumber
\|\mu_{n,\bom}-\mu_{\bom}\|&\leq&\|\tilde\mu_{n,p}\tau_{n,\bom}^{-1}-\tilde\mu_{p}\tau_{\bom}^{-1}\|\\
\nonumber
&\leq&\|\tilde\mu_{n,p}\tau_{n,\bom}^{-1}-\tilde\mu_{p}\tau_{n,\bom}^{-1}\|+\|\tilde\mu_{p}\tau_{n,\bom}^{-1}-\tilde\mu_{p}\tau_{\bom}^{-1}\|\\
\label{eq:rhs1}
&\leq&\|\tilde\mu_{n,p}-\tilde\mu_{p}\|+\|\tilde\mu_{p}\tau_{n,\bom}^{-1}-\tilde\mu_{p}\tau_{\bom}^{-1}\|.
\end{eqnarray}
From \eqref{eq:mu1}, the first term in \eqref{eq:rhs1} goes to $0$. 
\\

\noindent
For the second term in \eqref{eq:rhs1}, we study the convergence of $\Sigma_2^n(\bom)$ to $\Sigma_2^n(\bom)$. 
Recall that $a(h)$ is defined in \eqref{eq:limit} and that the characteristic function of $I$ is given in \eqref{eq:FCI}. 
Since 
\begin{eqnarray*}
\nonumber
|e^{iuy}-1-iuy\ind_{|y|\leq 1}|&\leq& 2\ind_{|y|>1}+(1/2) u^2y^2\ind_{|y|\leq 1}\\
\label{eq:L21}
&\leq& (2\ind_{|y|>1}+(1/2)u^2\ind_{|y|\leq 1})(|y|^2\wedge 1)
\end{eqnarray*}
we derive from \eqref{eq:FCI} that $\phi_{I-I_{n}}(u)\to 1$ for all fixed $u$. 
Thus, we have $I_{n}-I\stackrel{{\cal L}}{\lto} \delta_0$ and $I_{n}\stackrel{\pit}{\lto} I$. 
Next, since $T_i$ has a bounded density for all $i\geq 1$, 
condition $2$ entails $h_n(T_i,\Delta_i)^2\wedge 1 \to h(T_i,\Delta_i)^2\wedge 1$ in $L^1(\Omega, {\cal F}, \pit)$ and $h_n(T_i,\Delta_i)\stackrel{\pit}{\lto} h(T_i,\Delta_i)$. 
So that $\sum_{i=1}^ph_n(T_i,\Delta_i)\stackrel{\pit}{\lto} \sum_{i=1}^ph(T_i,\Delta_i)$. 
Together with $I_n\stackrel{\pit}{\lto}I$, we derive on the disintegrated probability space $\bar \Omega_{p+1}$:
$$
\Sigma_{p+1}^n:= \sum_{k\geq p+1} h_n(T_k,\Delta_k)\stackrel{\bar\pit}{\lto}\Sigma_{p+1}:= \sum_{k\geq p+1} h(T_k,\Delta_k).
$$ 
For any subsequence $(n')\subset (n)$, there is some further subsequence $(n'')\subset (n')$ and $\bar\Omega_{p+1}^0$ with $\bar\pit_{p+1}(\bar\Omega_{p+1}^0)=1$ 
such that for every $\bom\in\bar\Omega_{p+1}^0$, the convergence $\lim_{n''\to+\infty}\Sigma_{p+1}^{n''}(\bom)=\Sigma_{p+1}(\bom)$ holds true.
Since the operator of translation is continuous in $L^1(\rit)$, 
and since, under condition 3, $\tilde \mu_1\ll\lambda$, the second term in \eqref{eq:rhs1} goes also to $0$.  
This yields $\mu_{n'',\bom}\cvar \mu_{\bom}$ when $n''\to+\infty$ 
for all $\bom\in\bar\Omega_0$.
\\

\noindent 
Finally from the disintegration formula \eqref{eq:desint1}, for any $(n')\subset (n)$, there is some $(n'')\subset (n')$ such that  $\mu_{n''}\cvar \mu$. 
This proves $\mu_n\cvar \mu$. 
\CQFD
\end{Proof}

\begin{Rem}[Discussion on the conditions of Proposition \ref{prop:P4}]
\label{rem:var1} \

{\rm
\begin{itemize}

\item As seen before, if $h_n, h\in L^1(\lambda_{\rit_+}\otimes\sigma)$, then conditions \eqref{eq:Lévy2}--\eqref{eq:limit} are satisfied for the existence of $I_n$ and of $I$.
Moreover, if $\lim_{n\to+\infty}h_n=h$ in $L^1(\lambda_{\rit_+}\otimes\sigma)$, conditions~$1$ and $2$ hold true. 

\item Suppose $\sigma\ll\lambda^d$, 
then if, for almost all $t\in\rit$, $h_n(t,x)$ and $h(t,x)$ are absolutely continuous as functions of $x_i$ for all $1\leq i\leq d$ 
and $h_n\to h$, $\partial_{x_i}h_n(t,x) \to\partial_{x_i}h(t,x)$, $n\to+\infty$, in $L^r(\rit^{1+d})$ 
for some $r\geq 1+d$ with $\det( \partial_{x_i}h(t,x), \partial_{x_j}h(t,x))_{i,j=1}^d\neq 0$  
then Condition $3$  holds true in the simplest case where $p=1$. 
The same holds true if for almost all $x\in\rit^d$, $h_n(t,x)$ and $h(t,x)$ are absolutely continuous function of $t$
and $h_n\to h$, $\partial_th_n(t,x) \to\partial_th(t,x)$, $n\to+\infty$, in $L^r(\rit^{1+d})$ for some $r\geq 1+d$ 
with $\partial_th(t,x)\neq 0$ for almost all $(t,x)$, see Lemma \ref{prop:ABP}. 

%
%
\end{itemize}
 }
 \end{Rem}

\begin{Rem}{\bf (Comparison of Proposition \ref{prop:P4} with \cite{Kulik4})}
\label{rem:var1b} \

{\rm
Like in Remark \ref{rem:P1a1}, observe that the "time-stretching" transformation used by Kulik requires a smooth behavior of the kernels $h_n$ with respect to time. 
Indeed, the convergence of the stochastic derivates of $I_n$ assumed in \cite[Th. 2.2]{Kulik4} requires in particular 
the differentiability of the kernels $h_n(t,x)$ with respect to $t$ and a convergence of the derivatives. 
In contrast, Proposition \ref{prop:P4} above can be applied without smooth behavior of $h(t,x)$ with respect to $t$, 
see the second point in Remark \ref{rem:var1} above where explicit conditions in terms of smoothness of $h(t,x)$ but with respect to $x$ are given. 
 
}
\end{Rem}

 \bigskip
When $d=1$ and $h_n(t,x)=xg_n(t)$, we can adapt the proof of Proposition \ref{prop:P4} for the shot noise series $I_n=\sum_{k\geq 1} \Delta_kg_n(T_k)$ under more specific conditions.
This is in particular interesting when $h_n(t,x)=xg_n(t) \ind_{B_1^c}(x)$. 
In this case, the series $I_n$ becomes the Poisson integrals $\int_0^{+\infty}\int_{|x|>1} xg_n(t) N(dt,dx)$ appearing for instance in the Lévy-Itô decomposition of a Lévy-type stochastic integral \eqref{eq:gY}. 
The following result applies in particular in such setting: 
\begin{corollary}
\label{corol:P3}
For a law $\sigma$ with a finite first moment and $g_n \in L^1(\lambda_{\rit_+})$, 
consider the shot noise series $I_n=\sum_{k\geq 1} \Delta_kg_n(T_k)$. 
Suppose $g_n\to g$ in $L^1(\lambda_{\rit_+})$, $\lambda_{\rit_+}g^{-1}\ll\lambda$ and for all $t>0$
\begin{equation}
\label{eq:H3b}
\lambda_{[0,t]}g_n^{-1}\cvar \lambda_{[0,t]} g^{-1}. 
\end{equation}
Then $\disp {\cal L}(I_n)\cvar {\cal L}(I)$ when $n\to+\infty$.
\end{corollary}

\begin{Proof}
Under the conditions of Corollary \ref{corol:P3}, 
we have $xg_n(t)\to xg(t)$, $n\to+\infty$, in $L^1(\lambda_{\rit_+}\otimes\sigma)$ and Remark \ref{rem:var1} shows that the proof of Proposition \ref{prop:P4} still works with $p=1$ therein.
The only point to revise is \eqref{eq:mu1} (with $p=1$). 
But Condition \eqref{eq:H3b} entails ${\cal L}(g_n(U))\cvar {\cal L}(g(U))$ for all uniform random variable $U$, independent of $(\Delta_i)_{i\geq 1}$, 
and Lemma \ref{lemme:1} ensures
${\cal L}(\Delta_1g_n(U))\cvar {\cal L}(\Delta_1 g(U))$.
The rest of the proof follows the same lines as that of Proposition \ref{prop:P4} since ${\cal L}(\Delta_1 g(U))\ll\lambda$
when $\lambda_{\rit_+}g^{-1}\ll\lambda$ and $\sigma$ does not have atom in zero. 
\CQFD
\end{Proof}

\begin{Rem}
\label{rem:varg1}
{\rm 
If the functions $g_n, g$ are absolutely continuous with $g_n(0)\to g(0)$ and $g_n'\to g'$ in $L^1(\lambda_{\rit_+})$ when $n\to+\infty$ and with $g'(t)\neq 0$ for almost all $x\geq 0$, 
then the conditions \eqref{eq:H3b} and $\lambda_{\rit_+}g^{-1}\ll\lambda$ are satisfied, see Lemma \ref{lemme:Davydov3}.
}
\end{Rem}


\subsection{Convergence for truncated shot noise series}
\label{sec:varI(t)}
Truncated shot noise series $I(t)$ in \eqref{eq:series2} have an atom in zero and specific arguments have to be given to derive again   
the convergence in variation of the laws. 
We recall that this setting applies in particular to Poisson integrals in the Lévy-Itô decomposition of Lévy-type integrals on $[0,t]$
when $h_n(s,x)=f_n(s,x)\ind_{[0,t]\times B_1^c}(s,x)$, see \eqref{eq:gY}.
We use the following elementary result:
\begin{lemma}
\label{lemme:UnifStat}
Conditionally to $A_i=\{T_i\leq t< T_{i+1}\}$, the vector $(T_1,\dots, T_i)$ is the uniform order statistics, 
{\it i.e.} its law is given by the density $(i!/t^i)\ind_{0 \leq t_1\leq t_2\leq \cdots\leq t_i\leq t}$. 
\end{lemma}
The main result for the truncated shot noise series is: 
\begin{prop}
\label{prop:P3ter}
Let $f_n$ and $f$ 
be such that the shot noise series $I_n(t)$, $I(t)$ are well defined for some fixed $t>0$. 
Suppose 
\begin{equation}
\label{eq:H3bis}
(\lambda_{[0,t]}\otimes \sigma) f_n^{-1}\cvar (\lambda_{[0,t]}\otimes \sigma) f^{-1}.
\end{equation}
Then ${\cal L}(I_n(t))\cvar {\cal L}(I(t))$ when $n\to+\infty$. 
\end{prop}
\begin{Rem}
\label{rem:var2}
{\rm  \begin{itemize}
\item Like in Remark \ref{rem:var1}, explicit sufficient conditions for \eqref{eq:H3bis} are given by Lemma \ref{prop:ABP}. 
Moreover, explicit conditions for the existence of $I_n(t)$ are given in \eqref{eq:Lévy2} and \eqref{eq:limit} 
with $h(s,x)=f(s,x) \ind_{[0,t]}(s)$. 

\item In Proposition \ref{prop:P4}, the main condition \eqref{eq:cvar-main} asks, roughly speaking, for the convergence in variation of the law of some finite sub sum of the shot noise series. 
In the truncated case of Proposition \ref{prop:P3ter}, the series have only a finite (random) number of terms
 and possibly just one term with a positive probability. 
Consequently the condition \eqref{eq:cvar-main} does not seems sufficient, and instead we ask in \eqref{eq:H3bis} for the convergence in variation of the law of one term. 
\end{itemize}
}
\end{Rem} 

\begin{Proof}
The truncated series is made of a finite number of summands. 
We condition by the number of terms in the series and we deal with each case. 
Let $\mu_n={\cal L}(I_n(t))$ and $\mu={\cal L}(I(t))$ and for $A_0=\{T_1>t\}$ and $A_i=\{T_i\leq t<T_{i+1}\}$, $i\geq 1$, 
let $\mu_{n,i}={\cal L}(I_n(t)|A_i)$ and $\mu_i={\cal L}(I(t)|A_i)$. 
We have: 
\begin{equation*}
\label{eq:mumu}
\mu_n=\pit(A_0)\delta_0+\sum_{i=1}^{+\infty} \pit(A_i)\mu_{n,i}
\quad\mbox{ and } \quad
\mu=\pit(A_0)\delta_0+\sum_{i=1}^{+\infty} \pit(A_i)\mu_i
\end{equation*}
and for arbitrary $p\geq 1$:
\begin{eqnarray*}
\|\mu_n-\mu\|
&\leq& \sum_{i=1}^{+\infty} \pit(A_i) \|\mu_{n,i}-\mu_i\|
\leq \sum_{i=1}^{p} \pit(A_i) \|\mu_{n,i}-\mu_i\|+2\sum_{i=p+1}^{+\infty} \pit(A_i).
\end{eqnarray*}
Since $\sum_{i=0}^{+\infty} \pit(A_i)=1$ is a convergent series, 
it is enough to show, for all $i\geq 1$, $\mu_{n,i} \cvar \mu_i$ when $n\to +\infty$.
\\

\noindent
Since conditionally to $A_i$, $I_n(t)$ and $I(t)$ rewrites $I_n(t)=\sum_{k=1}^i f_n(T_k',\Delta_k)$ 
and $I(t)=\sum_{k=1}^i f(T_k',\Delta_k)$,  
Lemma \ref{lemme:UnifStat} (and commutativity of addition) entails that conditionally to $A_i$, $I_n(t)$ and $I(t)$ have the same law as
$\sum_{k=1}^i f_n(U_k,\Delta_k)$ and $\sum_{k=1}^i f(U_k,\Delta_k)$ 
where $U_k$ ($1\leq k\leq i$) are i.i.d. uniform random variables on $[0,t]$. 
By independence, the law of $\sum_{k=1}^i f_n(U_k,\Delta_k)$ is the convolution of the law of $f_n(U_k,\Delta_k)$, $1\leq k\leq i$, that is
\begin{equation*}
\label{eq:*in}
{\cal L}\left(\sum_{k=1}^i f_n(U_k,\Delta_k)\right)
=\big((t^{-1}\lambda_{[0,t]}\otimes\sigma)f_n^{-1} \big)^{*i}
\end{equation*}
and similarly
\begin{equation*}
\label{eq:*i}
{\cal L}\left(\sum_{k=1}^i f(U_k,\Delta_k)\right)
=\big((t^{-1}\lambda_{[0,t]}\otimes\sigma)f^{-1} \big)^{*i}.
\end{equation*}
Finally under condition \eqref{eq:H3bis}, Lemma \ref{lemme:convol} achieves the proof. 
\CQFD
\end{Proof}

\bigskip
In the case of integrands $f(s,x)=xg(s)$ related to Poisson integrals, Proposition \ref{prop:P3ter} rewrites as follows:
\begin{corollary}
\label{cor:2}
Let $g_n$ be such that for some fixed $t$, the shot noise series $I(t)$ with kernels $h_n(s,x):=xg_n(s)$ are well defined. 
Suppose that
\begin{equation}
\label{eq:H5}
\lambda_{[0,t]}g_n^{-1}\cvar\lambda_{[0,t]}g^{-1}.
\end{equation}
Then the laws ${\cal L}(I_n(t))$ converge in total variation to ${\cal L}(I(t))$.   
\end{corollary}
\begin{Rem}
\label{rem:varg2}
{\rm  
Sufficient conditions for \eqref{eq:H5} to hold true are given by Lemma \ref{lemme:Davydov3}.
}
\end{Rem}

\begin{Proof}
Like in the proof of Proposition \ref{prop:P3ter}, it is enough to show for all $i\geq 1$ that 
$$
{\cal L}\left(\sum_{k=1}^i \Delta_k g_n(U_k)\right)
\cvar 
{\cal L}\left(\sum_{k=1}^i \Delta_k g(U_k)\right), \med n\to+\infty.
$$
This comes from Lemma~\ref{lemme:1} and Lemma \ref{lemme:convol} when \eqref{eq:H5} holds true. 
\CQFD
\end{Proof}


\section{Application of the method to a stochastic differential equation}
\label{sec:EDS}

In this section, we are interested in the law of the solutions of  a Lévy driven stochastic differential equation (SDE)
\begin{equation}
\label{eq:EDS-X}
X_{n,t}=x_{n,0}+\int_0^t a_n(X_{n,s})ds +Z_t
\end{equation}
where $a_n$ are $C^1$ drifts with bounded derivatives and $Z$ is given in Section \ref{sec:shotnoise}. 
We apply our technique with a slight modification contrasting with the previous sections. 
Indeed, we apply a perturbation which consists in erasing the first jump of the driven Lévy process $Z$. 
However the technique is still elementary and deals only with the first jumps.  
We show that when the parameters of the SDE converge, we have the convergence of variation of the law of the solution. 
This recovers a particular case of Theorem 1.2 in \cite{Kulik3} (see also \cite[Th.~4.1]{Kulik4}). 

In such an SDE, Nourdin and Simon show in \cite{NS} the regularizing effect of the  drift term $a(X_s)ds$ on the law of the solution, {\it i.e.} 
when $a$ is locally monotonous at the initial condition $x_0$, 
${\cal L}(X_t)\ll \lambda \Longleftrightarrow {\cal L}(X_1)\ll \lambda \Longleftrightarrow \nu$ is infinite. 
Our initial motivation was to show that the regularizing effect of the drift actually works on the total variation on the law. 
This is the content of the main result of the section.
In a $n$-dimensional setting, when the drift is linear ($a(x)=a\cdot x$), 
the solution to \eqref{eq:EDS-X} is a Ornstein-Uhlenbeck process whose law has been studied in \cite{Simon} with similar technique as our (see also \cite{PZ}).   
\begin{prop}
\label{prop:P5}
Assume the Lévy measure $\nu$ is infinite and let $X_n$ be the solution of the SDE \eqref{eq:EDS-X} with drift function $a_n$ and initial condition $x_{n,0}$. 
Suppose 
\begin{enumerate}
\item $\lim_{n\to+\infty}x_{n,0}= x_0$.  
\item $a_n$ and $a$ are differentiable with bounded derivatives, 
the convergence $a_n'\to a'$ is uniform on bounded sets and for some fixed $t_0$, $a_n(t_0)\to a(t_0)$.
\item $a_n(y)$ and $a_n'(y)$ are both continuous functions of the couple $(n,y)$.
\item $a$ is locally monotonous at $x_0$. 
\end{enumerate}
Then, for each $t>0$, the law of $X_{n,t}$ converges in variation to that of $X_t$. 
\end{prop}

\begin{Proof} 
Since this result is proved in a greater generality in \cite{Kulik3}, we only sketch our elementary argument for $t=1$, following the setting in \cite{NS}.  
Since $\nu([0,1])=+\infty$, we can define a sequence of jumping times $T_i$ corresponding to small jumps $\Delta_i$ of $Z$
such that, with a probability close to $1$, $T_2<1$. 
We disintegrate the probability space as follows 
$\left(\bar\Omega_2\times[0, T_2(\bom)], \bar{\cal F}_2\times {\cal B}([0, T_2(\bom)]), \bar\pit_2\otimes \bar\lambda_{[0,T_2(\bom)]}\right)$
where $(\bar\Omega_2, \bar{\cal F}_2, \bar\pit_2)$ is the canonical space associated to $(\Delta_1, \bar Z)$ 
with $\bar Z_t=Z_t-\Delta_1\ind_{T_1\leq t}$ ($\bar Z$ is $Z$ with its first jump being erased). 
Consider $Y_n=X_n-Z$ and $Y=X-Z$  and observe that they are solutions of ordinary differential equations.  
For a subset $\bar\Omega_2^0 \in\bar{\cal F}_2$ of probability close to $1$, $X$ and $Y$ remain close to $x_0$ on $[0,T_2]$. 
Denoting $\mu_n$ and $\mu$ for the law of $X_{n,1}$ and of $X_1$, we have 
\begin{eqnarray*}
\nonumber 
\|\mu_n-\mu\|
\label{eq:var}
&\leq& 2\bar\pit_2((\bar\Omega_2^0)^c)+\int_{\bar \Omega_2^0}\|\bar\mu_{n,\bom}-\bar\mu_{\bom}\| d\bar\pit_2
\end{eqnarray*}
where 
$$
\bar\mu_{n,\bom}=
\bar{\lambda}_{[0,T_2(\bom)]}X_{n,1}(\bom, \cdot)^{-1}
\quad \mbox{ and } \quad
\bar\mu_{\bom}=
\bar{\lambda}_{[0,T_2(\bom)]}X_1(\bom, \cdot)^{-1}.
$$
We stress that on the disintegrated probability space $\bar\Omega_2\times [0,T_2(\bom)]$, the random variables $X_1$, $X_{n,1}$, $Y_1$, $Y_{n,1}$ and $Z_1$ 
are seen as variables of $(\bom, T_1)$. 
However, observe that $Z_1(\bom,\cdot)$ actually does not depend on $T_1$ since $Z$ jumps at least twice before $1$ 
and by the Lévy-Itô decomposition, the terminal value $Z_1$ is independent of the first jumping time $T_1$. 
By simplicity we shall write $Z_1(\bom)$.
Setting $\tau_{Z_1(\bom)}$ for the translation of $Z_1(\bom)$ in $\rit$, we rewrite 
\begin{eqnarray*}
\bar{\lambda}_{[0,T_2(\bom)]}X_{n,1}(\bom, \cdot)^{-1}
&=&\bar{\lambda}_{[0,T_2(\bom)]}Y_{n,1}(\bom, \cdot)^{-1} \tau_{Z_1(\bom)}^{-1}
\\
\bar{\lambda}_{[0,T_2(\bom)]}X_{1}(\bom, \cdot)^{-1}
&=&\bar{\lambda}_{[0,T_2(\bom)]}Y_{1}(\bom, \cdot)^{-1} \tau_{Z_1(\bom)}^{-1}
\end{eqnarray*}
and it remains to show for all $\bom\in\bar\Omega_2^0$:
\begin{equation*}
\label{eq:but}
\bar{\lambda}_{[0,T_2(\bom)]}Y_{n,1}(\bom, \cdot)^{-1}
\cvar
\bar{\lambda}_{[0,T_2(\bom)]} Y_1(\bom, \cdot)^{-1}.
\end{equation*} 
Note that under those notations when $\bom$ is fixed, $T_1$ is just a parameter in $[0,T_2(\bom)]$ 
and we are investigating for the convergence in variation of image measures by the mappings $Y_n(\bom, \cdot)$ on $[0, T_2(\bom)]$. 
This is proved by applying Lemma \ref{lemme:Davydov3} since from  \cite[Proposition 2]{NS}, $Y_{n,1}$ depends differentiably on the parameter $T_1$ with the derivative given by
\begin{equation*}
\label{eq:Y'}
\frac{dY_{n,1}}{dT_1}=
(a_n(X_{T_1^-})-a_n(X_{T_1})) \exp\left(\int_{T_1}^1 a_n'(X_{s}) ds\right).
\end{equation*}
Our hypothesis on $a_n$ and on $a_n'$ ensure that the first two conditions of Lemma \ref{lemme:Davydov3} are satisfied. 
The third one is ensured by the local monotonicity of $a$ at $x_0$ since $X_{T_1^-}$ and $X_{T_1}$ are closed to $x_0$ for $\bom\in\bar\Omega_2^0$. 
\CQFD
\end{Proof}
 

\bigskip
\noindent
{\bf Acknowledgment.} The author thanks anonymous referees for valuable comments 
which have permitted to propose a considerably improved version of this paper.

\end{document}